\documentclass[leqno,11pt]{article}
\usepackage{amssymb,latexsym,mathrsfs,color}
\usepackage[dvips]{graphicx} 
\usepackage{url}

 \catcode`@=11\def\@oddhead{\hbox{}\hfil\rm\thepage}\def\@oddfoot{}
 \def\@evenhead{\hbox{}\hfil\rm\thepage}\def\@evenfoot{}
 \@twosidetrue \catcode`@=12
\newtheorem{prp}{Proposition}
\newtheorem{lem}[prp]{Lemma}\newtheorem{thm}[prp]{Theorem}
\newtheorem{cor}[prp]{Corollary}
\newenvironment{prf}{\begin{trivlist}\item[\emph{Proof.}]}{\end{trivlist}
  \medskip\par}
\newenvironment{prfof}[1]{\begin{trivlist}\item[\emph{Proof of #1.}]}{
  \end{trivlist} \medskip \par}
\newenvironment{rem}{\begin{trivlist}\item[\emph{Remark.}]}{\end{trivlist}
  \medskip\par}
\def\prpb{\begin{prp}}\def\prpe{\end{prp}}
\def\lemb{\begin{lem}}\def\leme{\end{lem}}
\def\thmb{\begin{thm}}\def\thme{\end{thm}}
\def\corb{\begin{cor}}\def\core{\end{cor}}
\def\prfb{\begin{prf}}\def\prfe{\end{prf}}
\def\prfofb#1{\begin{prfof}{#1}}\def\prfofe{\end{prfof}}
\def\remb{\begin{rem}}\def\reme{\end{rem}}
\def\prpa#1{\label{p:#1}}\def\prpu#1{Proposition~\ref{p:#1}}
\def\lema#1{\label{l:#1}}\def\lemu#1{Lemma~\ref{l:#1}}
\def\thma#1{\label{t:#1}}\def\thmu#1{Theorem~\ref{t:#1}}
\def\cora#1{\label{c:#1}}\def\coru#1{Corollary~\ref{c:#1}}
\def\seca#1{\label{s:#1}}\def\secu#1{\S~\ref{s:#1}}

\def\itmb{\begin{enumerate}}\def\itme{\end{enumerate}}
\def\itdb{\begin{itemize}}\def\itde{\end{itemize}}
\def\ittb{\begin{description}}\def\itte{\end{description}}
\def\eqnb{\begin{equation}}\def\eqne{\end{equation}}
\def\arrb#1{\begin{array}{#1}}\def\arre{\end{array}}
\def\tabb#1{\par\noindent\begin{tabular}{#1}}
\def\tabe{\end{tabular}\par\noindent}
\def\eqna#1{\label{e:#1}}\def\eqnu#1{(\ref{e:#1})}

\def\QED{\relax\ifmmode\let\@tempa\relax\ifcase\@eqcnt\def\@tempa{& & &}\or
  \def\@tempa{& &}\else\def\@tempa{&}\fi\@tempa $\Box$ \else\hfill $\Box$ \fi}
\def\DDD{\relax\ifmmode\let\@tempa\relax\ifcase\@eqcnt\def\@tempa{& & &}\or
 \def\@tempa{& &}\else\def\@tempa{&}\fi\@tempa $\Diamond$
 \else\hfill $\Diamond$ \fi}

\def\Rom#1{\uppercase\expandafter{\romannumeral#1}}
\def\dsp{\displaystyle}

\def\limf#1{\displaystyle \lim_{#1\to\infty}}

\def\Ccomb#1#2{\setbox0=\hbox{$\displaystyle\mathrm{C}$}\setbox1=\hbox{%
$\scriptstyle #1$}\kern \wd1{\mathrm{C}}_{\kern -1.05\wd0\kern -0.99\wd1{#1}
 \kern 1.15\wd0{#2}}}

\def\clvec#1#2#3{\def\clvecone{#3}\left(\arrb{c} \dsp #1\\ \dsp #2
 \ifx\clvecone\empty\else\\ \dsp #3\fi\arre\right)}

\def\bar#1{\overline{#1}}
\def\le{\leqq} \def\ge{\geqq} 

\def\reals{{\mathbb R}}
\def\nreals#1{{\mathbb R}^{#1}}
\def\preals{[0,\infty)} 
\def\pintegers{{\mathbb Z}_+}
\def\nintegers{{\mathbb N}}

\def\prb#1{\def\prbone{#1}
  \ifx\prbone\empty{\mathrm{P}}\else{\mathrm{P[\;}}#1{\mathrm{\;]}}\fi}
\def\prbseq#1#2{\def\prbseqone{#2}
  \ifx\prbseqone\empty{\mathrm{P}}_{#1}\ignorespaces
  \else{\mathrm{P}}_{#1}{\mathrm{[\;}}#2{\mathrm{\;]}}\fi}
\def\EE#1{{\mathrm{E[\;}}#1{\mathrm{\;]}}}
\def\EEseq#1#2{\def\EEseqone{#2}
  \ifx\EEseqone\empty{\mathrm{E}}_{#1}\else
 {\mathrm{E}}_{#1}{\dsp\mathrm{[\;}}#2{\mathrm{\;]}}\fi}

\def\VVseq#1#2{\def\VVseqone{#2}
  \ifx\VVseqone\empty{\matrm{V}}_{#1}\else
 {\mathrm{V}}_{#1}{\dsp\mathrm{[\;}}#2{\mathrm{\;]}}\fi}

\def\ssN{^{(N)}}

\def\calD{\mathcal{D}}
\def\calF{\mathcal{F}}

\title{
Doubly uniform complete law of large numbers for independent point processes
}
\author{
Tetsuya Hattori
\footnote{ 
Supported by JSPS KAKENHI Grant Number 26400146,
and by Keio Gijuku Academic Development Funds.
}
\\ 
\small Laboratory of Mathematics, Faculty of Economics, Keio University, 
\\
\small Hiyoshi Campus, 4--1--1 Hiyoshi, Yokohama 223-8521, Japan
\\ \small URL: \url{http://web.econ.keio.ac.jp/staff/hattori/research.htm}
\\ \small email: \url{hattori@econ.keio.ac.jp}
} 
\date{\today}

\begin{document}

\setcounter{section}{0}

\setcounter{equation}{0}

\maketitle

\begin{abstract}
We prove a law of large numbers in terms of complete convergence of
independent random variables taking values in increments of monotone 
functions, with convergence uniform both in the initial and the final time.
The result holds also for the random variables taking values in 
functions of $2$ parameters which share similar monotonicity properties 
as the increments of monotone functions.
The assumptions for the main result are the H\"older continuity on 
the expectations as well as moment conditions, 
while the sample functions may contain jumps.
In particular, we can apply the results to 
point processes (counting processes) which lack Markov or martingale
type properties.
\end{abstract}


2000 MSC Primary 60F15; Secondary 60G55

\section{Introduction.}
\seca{1}

Let $T>0$, which we fix throughout this paper.
Put $\dsp \Delta=\{(t_1,t_2)\in\nreals2\mid 0\le t_1\le t_2\le T\}$,
and denote by $\calD$ the set of functions 
$\dsp z:\ \Delta\to\preals$ which, for $(t_1,t_2)\in\Delta$, 
is non-increasing in $t_1$ and non-decreasing in $t_2$\,,
and satisfies $\dsp z(t,t)=0$ for $t\in[0,T]$.
Let $(\Omega,\calF,\prb{})$ be a probability space.
In this paper, we consider $\calD$ valued random variables.
Here, we say that $Z$ is a $\calD$ valued
random variable, if, for each $(t_1,t_2)\in \Delta$,
$Z(t_1,t_2):\ \Omega\to\reals$ is almost surely 
a real valued Borel random variable.
An example of $\calD$ valued random variable is $Z$ defined by
$\dsp Z(t_1,t_2)=\tilde{Z}(t_2)-\tilde{Z}(t_1)$ for 
a point (counting) process $\tilde{Z}$.

In this paper we prove the following.
\thmb
\thma{main}
Let $r>0$ and $q>2$.
For each $N\in \nintegers$, 
let $Z\ssN_i$, $i=1,2,\ldots,N$,
be a sequence of independent, $\calD$ valued random variables,
and let $M\ssN_i$ be a positive real,
and $w\ssN_i$, $i=1,2,\ldots,N$, a non-negative sequence.
Assume the following for each $i=1,2,\ldots,N$ and $N\in\nintegers$:
\eqnb
\eqna{LLNmonotoneindepproc_Holderexp_ufmestim_with_param11}
\arrb{l}\dsp
\mbox{(i)}\ \ \EE{Z\ssN_i(0,T)^q}^{1/q}\le M\ssN,
\\ \dsp
\mbox{(ii)}\ \ |\EE{Z\ssN_i(t_1,t_2)-Z\ssN_i(s_1,s_2)}|
\le M\ssN w\ssN_i\,(|t_1-s_1|^r+|t_2-s_2|^r),
\\ \dsp \phantom{\mbox{(ii)}\ \ }
\ (s_1,s_2),\ (t_1,t_2)\in \Delta.
\arre
\eqne
Then the arithmetic average
$\dsp Y\ssN=\frac1N\sum_{i=1}^N Z\ssN_i$
satisfies
\eqnb
\eqna{LLNmonotoneindepproc_Holderexp_ufmestim_with_param12}
\arrb{l}\dsp
\EE{ \sup_{(t_1,t_2)\in \Delta}
 |Y\ssN(t_1,t_2)-\EE{Y\ssN(t_1,t_2)} |^{q}}
\\ \dsp {}
\le \frac{\dsp M\ssN{}^{q}2^{q-1}}{\dsp N^{q^2r/(2qr+2r+2)}}\,
 (C_q^q\,(2T\,\bar{w}\ssN{}^{1/r}+1)+2^{2q}),
\\ \dsp {}
N=N_0,N_0+1,\ldots,
\arre
\eqne
where we put
\eqnb
\eqna{LLNmonotoneindepproc_Holderexp_ufmestim_with_param_notation}
\bar{w}\ssN=\frac1N\sum_{i=1}^N w\ssN_i\,,
\eqne
and $N_0$ is the smallest integer
satisfying $N_0^{qr/(2qr+2r+2)}\ge 2$,
and $C_q$ is a positive constant depending only on $q$.

If in addition, $\{M\ssN,\,\bar{w}\ssN\}$ is bounded, 
and
\eqnb
\eqna{LLNmonotoneindepproc_with_2param_exponent}
(q^2-2q-2)r>2
\eqne
holds,
then 
\eqnb
\eqna{LLNmonotoneindepproc_with_param}
\limf{N} \sup_{(t_1,t_2)\in \Delta} \biggl|
\frac1N\sum_{i=1}^N (Z\ssN_i(t_1,t_2)-\EE{Z\ssN_i(t_1,t_2)})\biggr|
=0,\ a.e.,
\eqne
holds.
\DDD
\thme
We note that concerning the dependence of random variables 
for different $N$,
we are considering a law of large numbers in terms of complete convergence
of Hsu and Robbins \cite{HR47,Er49,Er50},
where we assume nothing about dependence or independence
between the variables with different $N$.
In contrast, the law of large numbers for random walks considers 
a special case where $Z\ssN_i=Z^{(i)}_i$, $N\ge i$.
Instead of the term \textit{strong law of large numbers},
which often refers to the case of random walks,
we shall in this paper refer to
\eqnu{LLNmonotoneindepproc_with_param} 
the \textit{complete law of large numbers}
(uniform in $2$ parameters $t_1$ and $t_2$).

Note also that we assume H\"older continuity properties 
\eqnu{LLNmonotoneindepproc_Holderexp_ufmestim_with_param11}(ii)
only on the expectation. In particular, the sample paths
may contain jumps.

As a straightforward application to processes $Z\ssN_i$ with
$1$ time variable, we can apply \thmu{main} to the increment
$\dsp Z\ssN_i(t_1,t_2)\mapsto Z\ssN_i(t_2)-Z\ssN_i(t_1)$ of
increasing processes (such as point processes) $Z\ssN_i$\,,
to obtain the following.
\corb
\cora{main}
Let $D_{\uparrow}=D_{\uparrow}([0,T],\reals)$ be the set of non-decreasing,
right continuous functions on a closed interval $[0,T]$ with left limits
$\dsp y(t-0):=\lim_{s\uparrow t} y(s)$.
Let $r>0$ and $q>2$.
For each $N\in \nintegers$, 
let $Z\ssN_i:\ \Omega\to D_{\uparrow}$, $i=1,2,\ldots,N$,
be a sequence of independent, $D_{\uparrow}$ valued random variables,
and let $M\ssN_i$ be a positive real,
and $w\ssN_i$, $i=1,2,\ldots,N$, a non-negative sequence.
Assume the following for each $i=1,2,\ldots,N$ and $N\in\nintegers$:
\eqnb
\eqna{LLNmonotoneindepproc_Holderexp_ufmestim_with_param11c}
\arrb{l}\dsp
\mbox{(i)}\ \ \EE{|Z\ssN_i(T)-Z\ssN_i(0)|^q}^{1/q}\le M\ssN,
\\ \dsp
\mbox{(ii)}\ \ |\EE{Z\ssN_i(t)}-\EE{Z\ssN_i(s)}|
\le M\ssN w\ssN_i |t-s|^r,
\ s,t\in [0,T],
\arre
\eqne
where $\bar{w}\ssN$ is as in
\eqnu{LLNmonotoneindepproc_Holderexp_ufmestim_with_param_notation}.

Then the arithmetic average
$\dsp Y\ssN(t)=\frac1N\sum_{i=1}^N Z\ssN_i(t)$
satisfies
\eqnb
\eqna{LLNmonotoneindepproc_Holderexp_ufmestim_with_param12c}
\arrb{l}\dsp
\EE{ \sup_{(t_1,t_2)\in \Delta}
 |Y\ssN(t_2)-Y\ssN(t_1)-\EE{Y\ssN(t_2)-Y\ssN(t_1)} |^{q}}
\\ \dsp {}
\le \frac{\dsp M\ssN{}^{q}2^{q-1}}{\dsp N^{q^2r/(2qr+2r+2)}}\,
 (C_q^q\,(2T\,\bar{w}\ssN{}^{1/r}+1)+2^{2q}),
\\ \dsp {}
N=N_0,N_0+1,\ldots,
\arre
\eqne
where $N_0$ is the smallest integer
satisfying $N_0^{qr/(2qr+2r+2)}\ge 2$,
and $C_q$ is a positive constant depending only on $q$.

If in addition, $\{M\ssN,\,\bar{w}\ssN\}$ is bounded, 
and \eqnu{LLNmonotoneindepproc_with_2param_exponent} holds,
then a doubly uniform complete law of large numbers
\eqnb
\eqna{LLNmonotoneindepproc_with_paramc}
\limf{N} \sup_{(t_1,t_2)\in \Delta} \biggl|
\frac1N\sum_{i=1}^N (Z\ssN_i(t_2)-Z\ssN_i(t_1))
-\EE{Z\ssN_i(t_2)-Z\ssN_i(t_1)})\biggr|
=0,\ a.e.,
\eqne
holds.
\DDD
\core
%
If $Z\ssN_i(t)-\EE{Z\ssN_i(t)}$, $i=1,2,\ldots,N$, are martingales in $t$, 
then we can apply Doob's inequalities to bound the supremums of
their sums in $t$
by the values at last time $t=T$, and we can directly apply
the law of large numbers for real valued random variables.
In contrast, 
we assume no structures with respect to $t$ other than monotonicity
for $Z\ssN_i(t)$ in \coru{main}.
We give a simple example of such situation in \secu{replacement}.

\bigskip\par\textbf{Acknowledgment.}
The author would like to thank Prof.\ M.~Takei for helpful discussions.

\section{Doubly uniform bounds for increasing processes.}
\seca{fSRPvarphiLLNdemonstration}

In this section we prove \thmu{main}.
\prpb
\prpa{monotonediff_finitegate2}
Assume that $m\in\calD$ satisfies $m(0,T)\le 1$ and
is globally H\"older continuous;
there exists positive constants $r$ and $C$ such that
\eqnb
\eqna{Holdercontinuity}
|m(t_1,t_2)-m(s_1,s_2)|\le C|t_1-s_1|^r+C|t_2-s_2|^r,
\ \ (s_1,s_2),\ (t_1,t_2)\in \Delta
\eqne
Then for any $n\in\nintegers$ there exists a finite set
$\dsp \Delta^*=\{ (t_{k,1},t_{k,2})\in \Delta\mid k=1,2,\ldots,K \}$,
satisfying 
\eqnb
\eqna{monotonediff_finitegate21}
 K\le (2n-1)\,T(C\,n)^{1/r}+1,
\eqne
such that 
\eqnb
\eqna{monotonediff_finitegate24}
\sup_{(t_1,t_2)\in \Delta} |y(t_1,t_2)-m(t_1,t_2)|
\le \bigvee_{k=1}^{K} |y(t_{k,1},t_{k,2})-m(t_{k,1},t_{k,2})|
+\frac2n
\eqne
holds for any $y\in \calD$,
where $\dsp \bigvee_{k=1}^{K}c_k$ denotes the largest value in $\{c_k\}$.
\DDD\prpe
The following lemma is the technical core of this section.
\lemb
\lema{monotonediff_finitegate2}
Assume that $m\in\calD$ satisfies $m(0,T)\le 1$ and
that there exists a positive constant $r$ and $C$ such that
\eqnu{Holdercontinuity} holds.
Then for any $n\in\nintegers$ there exists a finite set
$\dsp \Delta^*=\{ (t_{i,1},t_{i,2})\in \Delta\mid i=1,2,\ldots,K \}$,
of size $K$ satisfying a bound \eqnu{monotonediff_finitegate21},
such that for each $(t_1,t_2)\in \Delta$, 
there exists $(s_1,s_2)\in \Delta^*_n$ such that
\eqnb
\eqna{monotonediff_finitegate22}
\ s_1\le t_1\,,\ t_2\le s_2\,, \ \mbox{ and }
\ m(s_1,s_2)-\frac2n\le m(t_1,t_2)\,,
\eqne
and
\eqnb
\eqna{monotonediff_finitegate23}
\arrb{l}\dsp
\mbox{either }\ m(t_1,t_2)\le \frac1n\,,\ \mbox{ or, 
there exists $(u_1,u_2)\in \Delta^*_n$ such that }
\\ \dsp
\ t_1\le u_1\,,\ u_2\le t_2\,, \ \mbox{ and }
\ m(t_1,t_2)\le m(u_1,u_2)+\frac2n\,.
\arre
\eqne
\DDD
\leme
\prfb
Note that monotonicity in the definition of $\calD$
implies that the maximum value of $m$ is 
attained at $(0,T)$ and the minimum is $m(t,t)=0$, $t\in[0,T]$.
We omit the trivial case of $m$ being identically $0$,
and consider the case $m(0,T)>0$.

Fix $n$ and put
\eqnb
\eqna{monotonediff_finitegate2_prf2}
\Delta_{n,i}=\{ (t_1,t_2)\in \Delta
\mid \frac{i}n <m(t_1,t_2)\le \frac{i+1}n\},
\ \ i=0,1,2,\ldots,n-1.
\eqne
Since $m$ takes values in $[0,1]$,
there is a partition of $\Delta$ in terms of 
$\{\Delta_{n,i},\ i=0,\ldots,n-1\}$;
\eqnb
\eqna{monotonediff_finitegate2_prf3}
 \Delta=\{(t_1,t_2)\in \Delta
\mid m(t_1,t_2)=0\}\cup\bigcup_{i=0}^{n-1} \Delta_{n,i}\,.
\eqne
Denote by $\dsp i_{max}$ the index 
such that $\dsp (0,T)\in \Delta_{n,i_{max}}$.
Then since $m$ is continuous, the mean value theorem and monotonicity
imply
\eqnb
\eqna{monotonediff_finitegate2_prf6}
\arrb{l}\dsp
0\le i\le i_{max}\ \Leftrightarrow\ \Delta_{n,i}\ne\emptyset,
\\ \dsp
0\le m(t_1,t_2)\le \frac{i_{max}+1}n,\ (t_1,t_2)\in \Delta.
\arre
\eqne
Fix $i\in\{ 1,2,\ldots,i_{max}\}$, and define
$t_{k,1}$ and $t_{k,2}$ inductively in $k=1,2,\ldots$, by
\eqnb
\eqna{monotonediff_finitegate2_prf4}
\arrb{l}
t_{0,2}=T,
\\ \dsp
t_{k,1}
=\inf\{ s\in[0,T]\mid m(s,t_{k-1,2})\in \Delta_{n,i-1}\}
\\ \dsp \phantom{t_{k,1}}
=\inf\{ s\in[0,T]\mid m(s,t_{k-1,2})\le\frac{i}n\},\ k=1,2,\ldots,
\\ \dsp
t_{k,2}
=\inf\{ s\in[0,T]\mid m(t_{k,1},s)\in \Delta_{n,i-1}\}
\\ \dsp \phantom{t_{k,1}}
=\inf\{ s\in[0,T]\mid m(t_{k,1},s)>\frac{i-1}n\},k=1,2,\ldots.
\arre
\eqne
\begin{center}
\includegraphics[scale=0.85]{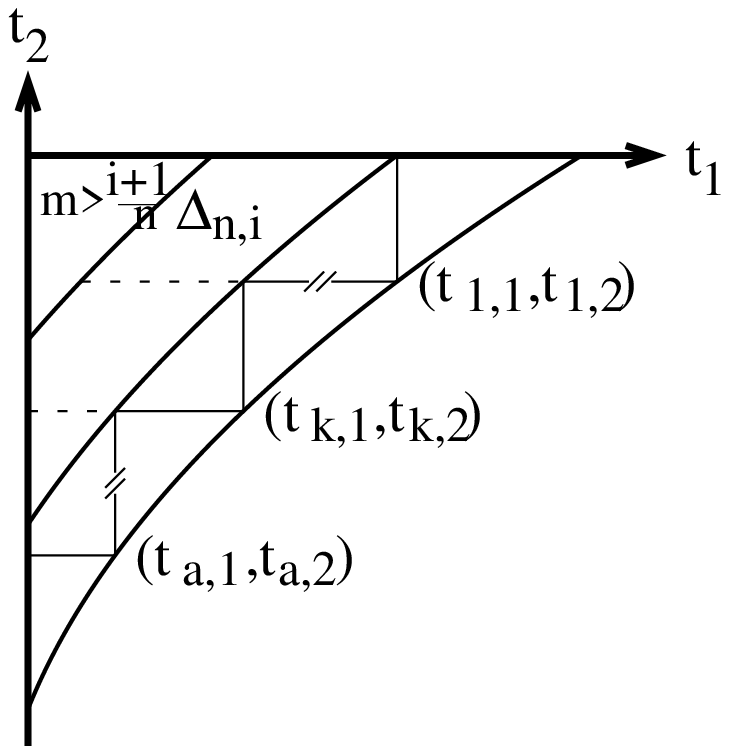}
\end{center}
Monotonicity of $m$ implies that 
$t_{k,1}$ and $t_{k,2}$ of \eqnu{monotonediff_finitegate2_prf4}
are respectively non-increasing in $k$,
and the continuity of $m$ implies
\eqnb
\eqna{monotonediff_finitegate2_prf7}
m(t_{k,1},t_{k,2})=\frac{i-1}n,
\ \ m(t_{k,1},t_{k-1,2})=\frac{i}n,
\ \ k=1,2,\ldots.
\eqne
This and the assumption of H\"older continuity \eqnu{Holdercontinuity} imply
\[
t_{k,2}-t_{k+1,2}\ge (C\,n)^{-1/r},\ k=0,1,2,\ldots,
\ \mbox{ and }\ 
t_{k,1}-t_{k+1,1}\ge (C\,n)^{-1/r},\ k=1,2,\ldots,
\]
so that if we denote the length of the sequence 
$(t_{k,1},t_{k,2})$, $k=1,2,\ldots$ by $a$,
then $\dsp a\le T\,(C\,n)^{1/r}$,
and the definition of $a$ also implies
$\dsp \frac{i-1}n\le m(0,t_{a,2})\le \frac{i}n$\,.

Let $(t_1,t_2)\in \Delta_{n,i}$\,.
Monotonicity, $\dsp m(0,t_{a,2})\le \frac{i}n$, and
\eqnu{monotonediff_finitegate2_prf2} imply
\eqnb
\eqna{monotonediff_finitegate2_prf8}
\exists k\in \{1,2,\ldots,a\};\ t_{k,2}< t_2\le t_{k-1,2}\,.
\eqne
This and monotonicity and $\dsp (t_1,t_2)\in \Delta_{n,i}$ and 
\eqnu{monotonediff_finitegate2_prf7} together imply
\[
m(t_{k,1},t_2)\le m(t_{k,1}, t_{k-1,2})=\frac{i}n<m(t_1,t_2),
\]
which further implies
\eqnb
\eqna{monotonediff_finitegate2_prf9}
 t_1<t_{k,1}
\eqne
and 
\eqnb
\eqna{monotonediff_finitegate2_prf5}
m(t_1,t_2)\le\frac{i}n+\frac1n=m(t_{k,1},t_{k,2})+\frac2n\,.
\eqne

Denote by
$\dsp \Delta^*_{+}=\{ (t_{\ell,1},t_{\ell,2})\in \Delta
\mid \ell=1,2,\ldots, a'\}$,
the union of so obtained 
$\{(t_{k,1},t_{k,2})\mid k=1,\ldots,a\}$
for all $i=1,\ldots,i_{max}$.
In particular, since $\dsp a\le T\,(C\,n)^{1/r}$ for each $i$,
$\dsp a'\le n\, T\,(C\,n)^{1/r}$.
Combining
\eqnu{monotonediff_finitegate2_prf3},
\eqnu{monotonediff_finitegate2_prf6},
\eqnu{monotonediff_finitegate2_prf8},
\eqnu{monotonediff_finitegate2_prf9}, and
\eqnu{monotonediff_finitegate2_prf5}, we have
\eqnb
\eqna{monotonediff_finitegate2_prf1}
\arrb{l}\dsp
(\forall (t_1,t_2)\in \Delta\setminus(\{(t_1,t_2)\in \Delta
\mid m(t_1,t_2)=0\}
\cup \Delta_{n,0})\,)
\\ \dsp
\ \exists (u_1,u_2)\in \Delta^*_{+};\ t_1< u_1,\ t_2\ge u_2\ \mbox{ and }\ 
m(t_1,t_2)\le m(u_1,u_2)+\frac2n\,.
\arre
\eqne
This completes a proof of \eqnu{monotonediff_finitegate23}.

A proof of the remaining bound is similarly proved.
Fix $n$, and put
\eqnb
\eqna{monotonediff_finitegate2_prf11}
\Delta'_{n,i}=\{ (t_1,t_2)\in \Delta \mid
 \frac{i}n \le m(t_1,t_2)<\frac{i+1}n\},\ \ i=0,1,\ldots,n.
\eqne
Then 
\eqnb
\eqna{monotonediff_finitegate2_prf12}
 \Delta=\bigcup_{i=0}^{n} \Delta'_{n,i}\,.
\eqne
Define an integer $i_{max}\le n$ by
\eqnb
\eqna{monotonediff_finitegate2_prf13}
0\le i\le i_{max}\ \Leftrightarrow\ \Delta'_{n,i}\ne\emptyset.
\eqne
Then
\eqnb
\eqna{monotonediff_finitegate2_prf14}
0\le m(t_1,t_2)< \frac{i_{max}+1}n\,,\ (t_1,t_2)\in \Delta.
\eqne
Fix $i\in\{ 0,1,2,\ldots,i_{max}-2\}$, and define
$t'_{k,1}$ and $t'_{k,2}$ inductively in $k=1,2,\ldots$, by
\eqnb
\eqna{monotonediff_finitegate2_prf15}
\arrb{l}
t'_{1,2}=T,
\\ \dsp
t'_{k,1}
=\inf\{ s\in[0,T]\mid m(s,t'_{k,2})\in \Delta'_{n,i+1}\}
\\ \dsp \phantom{t_{k,1}}
=\inf\{ s\in[0,T]\mid m(s,t'_{k,2})<\frac{i+2}n\},\ k=1,2,\ldots,
\\ \dsp
t'_{k,2}
=\inf\{ s\in[0,T]\mid m(t'_{k-1,1},s)\in \Delta'_{n,i+1}\}
\\ \dsp \phantom{t_{k,1}}
=\inf\{ s\in[0,T]\mid m(t'_{k-1,1},s)\ge \frac{i+1}n\},\ k=2,3,\ldots.
\arre
\eqne
\begin{center}
\includegraphics[scale=0.85]{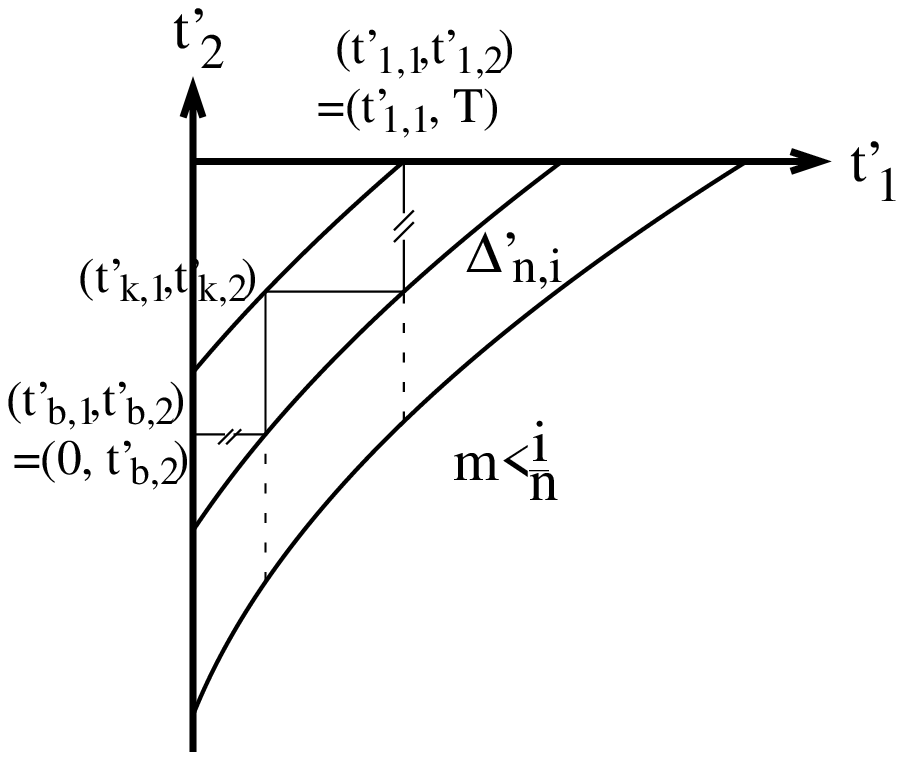}
\end{center}
Monotonicity of $m$ implies that 
$t'_{k,1}$ and $t'_{k,2}$ are respectively non-increasing in $k$,
and the continuity of $m$ implies
\eqnb
\eqna{monotonediff_finitegate2_prf16}
m(t'_{k,1},t'_{k,2})=\frac{i+2}n,
\ \ m(t'_{k,1},t'_{k+1,2})=\frac{i+1}n,
\ \ k=1,2,\ldots.
\eqne
This and the assumption of H\"older continuity \eqnu{Holdercontinuity} 
imply,
as before,
that there exists an integer $b$ satisfying
$\dsp b\le T\,(C\,n)^{1/r}$ and
$\dsp \frac{i+1}n\le m(0,t'_{b,2})\le \frac{i+2}n$\,,
such that the sequences $t'_{k,1}$ and $t'_{k,2}$ in $k$ are
of length $b$.
Put $t'_{b,1}=0$ for notational convenience below.

Let $(t_1,t_2)\in \Delta'_{n,i}$\,. Then 
\eqnb
\eqna{monotonediff_finitegate2_prf17}
t_1\ge t'_{1,1}\ \mbox{ or }
\exists k\in\{2,3,\ldots,b\};\ t'_{k,1}\le t_1< t'_{k-1,1}\,.
\eqne
This and monotonicity and $(t_1,t_2)\in \Delta'_{n,i}$ and
\eqnu{monotonediff_finitegate2_prf16} together imply
\[
m(t_1,t_2)\le \frac{i+1}n=m(t'_{k-1,1}, t'_{k,2})< m(t_1,t'_{k,2}),
\]
which further implies
\eqnb
\eqna{monotonediff_finitegate2_prf18}
 t_2<t'_{k,2}
\eqne
and
\eqnb
\eqna{monotonediff_finitegate2_prf19}
m(t'_{k,1},t'_{k,2}) = \frac{i+2}n\le m(t_1,t_2)+\frac2n\,.
\eqne
Denote
by $\dsp \Delta^*_{-}=\{ (t'_{\ell,1},t'_{\ell,2})\in \Delta
\mid \ell=1,2,\ldots, b'\}$,
the union of $\{(0,T)\}$ and so obtained 
$\{(t'_{k,1},t'_{k,2})\mid k=1,\ldots,b\}$
for all $i=0,1,\ldots,i_{max}-2$.
In particular, since $\dsp b\le T\,(C\,n)^{1/r}$ for each $i$,
and $i_{max}-2\le n-2$, we have
$\dsp b'\le (n-1)\, T\,(C\,n)^{1/r}+1$.
Combining
\eqnu{monotonediff_finitegate2_prf12},
\eqnu{monotonediff_finitegate2_prf13},
\eqnu{monotonediff_finitegate2_prf14},
\eqnu{monotonediff_finitegate2_prf17},
\eqnu{monotonediff_finitegate2_prf18}, and
\eqnu{monotonediff_finitegate2_prf19}, we have
\eqnb
\eqna{monotonediff_finitegate2_prf20}
\arrb{l}\dsp
(\forall (t_1,t_2)
\in \Delta\setminus(\Delta'_{n,i_{max}-1}\cup \Delta'_{n,i_{max}})
\ \exists (s_1,s_2)\in \Delta^*_{-};
\\ \dsp
t_1\ge s_1,\ t_2< s_2\ \mbox{ and }\ m(t_1,t_2)\ge m(s_1,s_2)-\frac2n\,.
\arre
\eqne
If $\dsp(t_1,t_2)\in \Delta'_{n,i_{max}}\cup \Delta'_{n,i_{max}-1}$,
then $t_1\ge 0$ and $t_2\le T$, and also
\eqnu{monotonediff_finitegate2_prf11},
\eqnu{monotonediff_finitegate2_prf13}, and the monotonicity of $m$
imply
$\dsp m(t_1,t_2)\ge m(0,T)-\frac2n$\,.
Since we included $(0,T)$ explicitly in $\Delta^*_-$, 
the existence of $(s_1,s_2)$ satisfying the estimates 
in \eqnu{monotonediff_finitegate2_prf20} holds also for this case.
This completes a proof of $\dsp m(s_1,s_2)-\frac2n\le m(t_1,t_2)$ in
\eqnu{monotonediff_finitegate22}.

Finally, we put $\Delta^*_n=\Delta^*_{+}\cup \Delta^*_{-}$. Then
\[
K=\sharp \Delta^*_n\le a'+b'\le (2n-1)(T\,(C\,n)^{1/r}+1)+1,
\]
which proves \eqnu{monotonediff_finitegate21}.
\QED\prfe
\prfofb{\protect\prpu{monotonediff_finitegate2}}
Let $(t_1,t_2)\in \Delta$.
If $y(t_1,t_2)\ge m(t_1,t_2)$, then,
\lemu{monotonediff_finitegate2} and monotonicity of $y$
imply that, there exists $(s_1,s_2)\in \Delta^*_n$ such that,
\[ \arrb{l}\dsp
|y(t_1,t_2)-m(t_1,t_2)|=y(t_1,t_2)-m(t_1,t_2)
\\ \dsp {}
\le (y(s_1,s_2)-m(s_1,s_2))+(m(s_1,s_2)-m(t_1,t_2))
\\ \dsp {}
\le |y(s_1,s_2)-m(s_1,s_2)|+\frac2n\,.
\arre\]
Assume in the following that $y(t_1,t_2)\le m(t_1,t_2)$.
Note that the assumptions on $y$ implies that $y$ is non-negative.
If $\dsp m(t_1,t_2)\le \frac1n$\,, then
\lemu{monotonediff_finitegate2} and nonnegativity of $y$ imply that
\[
|y(t_1,t_2)-m(t_1,t_2)|=m(t_1,t_2)-y(t_1,t_2)\le \frac1n\,.
\]
If $\dsp m(t_1,t_2)> \frac1n$\,, then
\lemu{monotonediff_finitegate2} and monotonicity of $y$
imply that, there exists $(u_1,u_2)$ in $\Delta^*_n$ such that,
\[ \arrb{l}\dsp
|y(t_1,t_2)-m(t_1,t_2)|
\\ \dsp {}
\le (m(u_1,u_2)-y(u_1,u_2))+(m(t_1,t_2)-m(u_1,u_2))
\\ \dsp {}
\le |y(u_1,u_2)-m(u_1,u_2)|+\frac2n\,.
\arre\]
Therefore \eqnu{monotonediff_finitegate24} holds.

\lemu{monotonediff_finitegate2} also implies the claimed bound on $K$,
the size of $\Delta^*_n$\,.
\QED
\prfofe
\corb
\cora{monotoneprocufmestimate2}
Let $Y$ be a random variable taking values in $\calD$.
Assume that
\eqnb
\eqna{Holdercontinuity1}
\arrb{l}\dsp
|\EE{Y(t_1,t_2)-Y(s_1,s_2)}|\le C\EE{Y(0,T)}\,(|t_1-s_1|^r+|t_2-s_2|^r),
\\ \dsp
(s_1,s_2),\ (t_1,t_2)\in \Delta,
\arre
\eqne
holds for some $r>0$ and $C>0$.
Then for any $n\in\nintegers$ and for any $q\ge1$,
\eqnb
\eqna{monotonediscretizationbd2}
\arrb{l}\dsp
\EE{\sup_{(t_1,t_2)\in \Delta} |Y(t_1,t_2)-\EE{Y(t_1,t_2)}|^{q}}
\\ \dsp {}
\le 2^{q-1} ((2n-1)\,T\,(C\,n)^{1/r}+1)
\sup_{(t_1,t_2)\in \Delta} \EE{|Y(t_1,t_2)-\EE{Y(t_1,t_2)}|^{q}}
\\ \dsp \phantom{\le}
+2^{q-1}\biggl(\frac2n\biggr)^{q}\EE{Y(0,T)}^q
\arre
\eqne
holds.
\DDD
\core
\prfb
For each sample $\omega\in\Omega$,
\prpu{monotonediff_finitegate2} with
\[
 y(t_1,t_2)=\frac{Y(\omega)(t_1,t_2)}{\EE{Y(0,T)}}
\,,\ 
m(t_1,t_2)=\frac{\EE{Y(t_1,t_2)}}{\EE{Y(0,T)}}
\]
implies that for any $n\in\nintegers$ there exists $K\ge1$ and 
a finite set
\[
\Delta^*=\{ (t_{k,1},t_{k,2})\in \Delta\mid k=1,2,\ldots,K \},
\]
independent of sample $\omega$, such that
\eqnb
\eqna{monotoneprocufmestimate2prf}
 K\le (2n-1)\,T(C\,n)^{1/r}+1,
\eqne
and
\[\arrb{l}\dsp
\sup_{(t_1,t_2)\in \Delta}
 |Y(\omega)(t_1,t_2)-\EE{Y(t_1,t_2)}|
\\ \dsp {}
\le \bigvee_{k=1}^{K} |Y(\omega)(t_{k,1},t_{k,2})-\EE{Y(t_{k,1},t_{k,2})}|
+\frac2n\,\EE{Y(0,T)}.
\arre\]
This, with an elementary inequality
\eqnb
\eqna{integopnonlinN}
a^p+b^p\le (a+b)^p\le 2^{p-1}(a^p+b^p) ,
\eqne
valid for all positive $a$ and $b$ with $p\ge1$, implies
\[\arrb{l}\dsp
\EE{\sup_{(t_1,t_2)\in \Delta} |Y(t_1,t_2)-\EE{Y(t_1,t_2)}|^q}
\\ \dsp {}
\le 2^{q-1}
\EE{\biggl( \bigvee_{k=1}^{K}
 |Y(t_{k,1},t_{k,2})-\EE{Y(t_{k,1},t_{k,2})}|
\biggr)^q}
+2^{q-1}\biggl(\frac2n\biggr)^q\EE{Y(0,T)}^q
\\ \dsp {}
= 2^{q-1} 
\EE{ \bigvee_{k=1}^{K} |Y(t_{k,1},t_{k,2})-\EE{Y(t_{k,1},t_{k,2})}|^q }
+2^{q-1}\biggl(\frac2n\biggr)^q\EE{Y(0,T)}^q
\\ \dsp {}
\le 2^{q-1} \sum_{k=1}^{K}
\EE{ |Y(t_{k,1},t_{k,2})-\EE{Y(t_{k,1},t_{k,2})}|^q }
+2^{q-1}\biggl(\frac2n\biggr)^q\EE{Y(0,T)}^q,
\arre\]
which, with \eqnu{monotoneprocufmestimate2prf}, implies
\eqnu{monotonediscretizationbd2}.
\QED
\prfe

We will see that the results so far reduce the claim of \thmu{main}
to the complete law of large numbers for real valued random variables.
It is known \cite{HR47,Er49,Er50} that 
a necessary and sufficient condition for the complete convergence of 
identically distributed real valued random variables with finite expectation
is the existence of variance.
Stated in terms of a complete law of large numbers,
existence of variance or the second order moment suffices for
the strong law to hold \cite[\S10.4, Example 1]{ChowTeicher}.
A generalized result for the case of varying distribution is also known
\cite{TH87}. We will use the results in the following form.
See the references for a proof.
\prpb
\prpa{LLNindepbdd}
For each $N\in\nintegers$, 
let $\tilde{Z}\ssN_i:\ \Omega\to\reals$, $i=1,2,\ldots,N$, 
be a finite sequence
of independent, real valued random variables, and put
$\dsp \tilde{Y}\ssN=\frac1N \sum_{i=1}^N \tilde{Z}\ssN_i$.
Assume that there exists $q>2$ such that
\[
M\ssN:= 
\max_{i\in \{1,\ldots,N\}} \EE{|\tilde{Z}\ssN_i|^q}^{1/q}<\infty.
\]
Then there exists a positive constant $C_q$ depending only on $q$,
(in particular, independent of $N$ and $M\ssN$,)
such that
\eqnb
\eqna{LLNindepbdd}
\EE{|\tilde{Y}\ssN-\EE{\tilde{Y}\ssN}|^q}^{1/q}
\le \frac{C_q M\ssN}{\sqrt{N}},
\eqne
hold.
\DDD
\prpe
We can for example put
\[
C_q=\biggl(\frac12(4k)^q+\frac{2k}{2k-q}(8k)^q\biggr)^{1/q},
\]
in \eqnu{LLNindepbdd}, where 
$k$ is the smallest integer greater than $q/2$.

We are ready to prove \thmu{main},
the complete law of large numbers for the arithmetic average
of $Z\ssN_i(t_1,t_2)$, uniform both in $t_1$ and $t_2$\,.
\prfofb{\protect\thmu{main}}
Note first that the assumption
\eqnu{LLNmonotoneindepproc_Holderexp_ufmestim_with_param11}(i)
implies
\eqnb
\eqna{LLNmonotoneindepproc_Holderexp_ufmestim_with_param_rem10}
\sup_{\stackrel{\scriptstyle i=1,2,\ldots,N,}{(t_1,t_2)\in \Delta}}
\EE{|Z\ssN_i(t_1,t_2)|^q}^{1/q}\le M\ssN,
\eqne
because of monotonicity.

The assumption
\eqnu{LLNmonotoneindepproc_Holderexp_ufmestim_with_param11}(ii)
imply that
\[
Y:=\frac1{M\ssN}Y\ssN=\frac1{N\,M\ssN}\sum_{i=1}^N Z\ssN_i
\]
satisfies all the assumptions in
\coru{monotoneprocufmestimate2},
with $\dsp C=\bar{w}\ssN$ in \eqnu{Holdercontinuity1}.
\coru{monotoneprocufmestimate2} and \prpu{LLNindepbdd} therefore imply,
\[\arrb{l}\dsp
\EE{ \sup_{(t_1,t_2)\in \Delta}
 |Y\ssN(t_1,t_2)-\EE{Y\ssN(t_1,t_2)} |^{q}}
\\ \dsp {}
\le 2^{q-1}((2n-1)T(\bar{w}\ssN\,n)^{1/r}+1)
\,\frac{(C_q \,M\ssN)^q}{N^{q/2}}
+2^{q-1}\biggl(\frac{2M\ssN}n\biggr)^{q}
\arre\]
for positive integers $n$ and $N$.
Now for each $N$, fix $n=n_N$ to be the largest integer not greater than 
$\dsp N^{rq/(2qr+2r+2)}$.
If $N_0^{rq/(2qr+2r+2)}\ge 2$, then for $N\ge N_0$,
 $\dsp \frac12N^{rq/(2qr+2r+2)}< n_N \le N^{rq/(2qr+2r+2)}$,
and we have
\eqnu{LLNmonotoneindepproc_Holderexp_ufmestim_with_param12}.

The assumption \eqnu{LLNmonotoneindepproc_with_2param_exponent}
and just proved result
\eqnu{LLNmonotoneindepproc_Holderexp_ufmestim_with_param12}
imply 
\[
\EE{ \sum_{N=N_0}^{\infty} \sup_{(t_1,t_2)\in \Delta}
 |Y\ssN(t_1,t_2)-\EE{Y\ssN(t_1,t_2)} |^{q}} <\infty,
\]
which implies \eqnu{LLNmonotoneindepproc_with_param}.
\QED
\prfofe

If the moment condition 
\eqnu{LLNmonotoneindepproc_Holderexp_ufmestim_with_param11}(i) holds for 
arbitrarily large exponent $q$,
the doubly uniform complete law of large numbers holds
with `order of fluctuation' arbitrary close to $1/2$, as expected.
\thmb
\thma{main2}
For each $N\in \nintegers$, 
let $Z\ssN_i:\ \Omega\to \calD$, $i=1,2,\ldots,N$,
be a sequence of independent, $\calD$ valued random variables.
Let $r>0$, and 
for $N\in\nintegers$, let $M\ssN$ be a positive real
and $w\ssN_i$, $i=1,2,\ldots,N$, a non-negative sequence.
Assume the following for each $i=1,2,\ldots,N$ and $N\in\nintegers$:
\eqnb
\eqna{LLNmonotoneindepproc_Holderexp_ufmestim_with_param1}
\arrb{l}\dsp
\mbox{(i)}\ \ \EE{|Z\ssN_i(0,T)|^q}^{1/q}\le M\ssN,
\ q\in\nintegers,
\\ \dsp
\mbox{(ii)}\ \ |\EE{Z\ssN_i(t_1,t_2)}-\EE{Z\ssN_i(s_1,s_2)}|
\le M\ssN w\ssN_i\, (|t_1-s_1|^r+|t_2-s_2|^r),
\\ \dsp \phantom{\mbox{(ii)\ \ }}
\ (s_1,s_2),\ (t_1,t_2)\in \Delta.
\arre
\eqne
Then for any $\dsp \gamma\in(0,\frac12)$, and $p>0$
the average 
$\dsp Y\ssN=\frac1N\sum_{i=1}^N Z\ssN_i$ satisfies
\eqnb
\eqna{main2}
\arrb{l}\dsp
\EE{ \sup_{(t_1,t_2)\in \Delta} 
|Y\ssN(t_1,t_2)-\EE{Y\ssN(t_1,t_2)} |^{p}}^{1/p}
\\ \dsp {}
\le \frac{M\ssN}{N^{\gamma}}2^{1-1/q}
 (C_q^q\,(2T(\bar{w}\ssN)^{1/r}+1)+2^{2q})^{1/q},
\\ \dsp {}
N=N_0,N_0+1,\ldots,
\arre
\eqne
where $\bar{w}\ssN$ is as in
\eqnu{LLNmonotoneindepproc_Holderexp_ufmestim_with_param_notation},
$C_q$ as in \prpu{LLNindepbdd},
$N_0=N_0(r,q)$ is the smallest integer
satisfying $N_0^{rq/(2rq+2r+2)}\ge 2$, and
$\dsp q=q(p,\gamma)= 3\vee \frac{r+1}r \,\frac{2\gamma}{1-2\gamma} \vee p$.
(In particular, $q$ and $N_0$ are independent of 
$N$, $M\ssN$, and $\bar{w}\ssN$.)

If in addition, $\{M\ssN,\,\bar{w}\ssN\}$ is bounded, 
then a strong uniform law of large numbers
\eqnu{LLNmonotoneindepproc_with_param}
holds.
\DDD
\thme
\prfb
As in the proof of \thmu{main},
we have \eqnu{LLNmonotoneindepproc_with_param} for all $N$ and $q$.

Let $\dsp \gamma\in(0,\frac12)$ and $p>0$,
and choose 
$\dsp q=3\vee 2\frac{r+1}r \,\frac{\gamma}{1-2\gamma} \vee p$.
Then $\dsp q\ge 2\frac{r+1}r \,\frac{\gamma}{1-2\gamma}$ implies
$\dsp \frac{rq}{2rq+2r+2}\ge \gamma$,
hence,
the monotonicity of $L_p$ norms and
\eqnu{LLNmonotoneindepproc_with_param}
imply
\[ \arrb{l}\dsp
\EE{ \sup_{(t_1,t_2)
\in \Delta} |Y\ssN(t_1,t_2)-\EE{Y\ssN(t_1,t_2)} |^{p}}^{1/p}
\\ \dsp {}
\le 
\EE{ \sup_{(t_1,t_2)\in \Delta}
 |Y\ssN(t_1,t_2)-\EE{Y\ssN(t_1,t_2)} |^{q}}^{1/q}
\\ \dsp {}
\le \frac{M\ssN}{N^{\gamma}}2^{1-1/q}
 (C_q^q\,(2T(\bar{w}\ssN)^{1/r}+1)+2^{2q})^{1/q},
\arre\]
which proves \eqnu{main2}.
The argument for the proof of \eqnu{LLNmonotoneindepproc_with_param}
in the proof of \thmu{main} also proves 
\eqnu{LLNmonotoneindepproc_with_param} in \thmu{main2}.
\QED
\prfe

\section{Example of point processes with dependent increments.}
\seca{replacement}

Results without assuming extra structures for counting processes
as in \coru{main} could be of use in practical situations.  For example,
consider a new large office building with a large number, say $N$, 
of lighting equipments.
Each light bulb has a random lifetime, 
which may depend on the location ($i=1,2,\ldots,N$) in the building.
The distribution of each lifetime also could depend in a mathematically
cumbersome way on the latest time the light bulb burnt out, 
because the light bulb products in the market are updated according to
e.g., advances in technology or regulation on materials.
We would be interested in estimating the number of bulbs to be 
replaced in the time period $[t_1,t_2]$, which is the random number
$N\, (Y\ssN(t_2)-Y\ssN(t_1))$ in terms of the notations in \coru{main}.
(Note also that in these practical applications where $N$ is finite
and fixed, a complete law of large numbers as we consider should be
natural than a strong law of large numbers which assumes relation
between random variables with different $N$.)

As an example of point process with dependent increments
we consider the point process $Z\ssN_i(t)$
with last-arrival-time dependent intensity \cite[\S3]{fluid14}.
For each $N\in\nintegers$ and $i=1,2,\ldots,N$, the
sequence $\tau\ssN_{i,0}=0<\tau\ssN_{i,1}<\cdots$ of arrival times
are random times defined inductively by
\begin{equation}
\eqna{intensity}
\begin{array}{l}\dsp
\prb{t< \tau\ssN_{i,k}\mid \calF_{\tau\ssN_{i,k-1}}}
=\exp(-\int_{\tau\ssN_{i,k-1}}^t w\ssN_i(\tau\ssN_{i,k-1},u)\,du)
\ \mbox{ on }\ t\ge \tau\ssN_{i,k-1}\,,
\end{array}
\end{equation}
where, $w\ssN_i$ is a non-negative continuously
differentiable function defined on
\[
 (t_0,t)\in \Delta=\{(t_1,t_2)\in\nreals2\mid 0\le t_1\le t_2\le T\},
\]
to which we refer as the intensity function of the counting process 
\eqnb
\eqna{PPLATDI}
Z\ssN_i(t)=\max\{k\in\pintegers\mid \tau\ssN_{i,k}\le t\}.
\eqne
If $w\ssN_i$ is independent of the first variable, the definition of
$Z\ssN_i$ reduces to that of the Poisson process with intensity function 
$w\ssN_i$, but in general, unlike the Poisson processes, 
$Z\ssN_i$ is not of independent increment.
\thmb
\thma{PPLATDILLN}
For $N\in\nintegers$ and $i=1,\ldots,N$, let $w\ssN_i:\ \Delta\to\preals$ be
a non-negative continuously differentiable function,
and $Z\ssN_i$ a process determined by \eqnu{PPLATDI}.
If $Z\ssN_i$, $i=1,\ldots,N$, are independent for each $N\in\nintegers$, 
and 
$\dsp C:= \sup_{N\in\nintegers,\ i\in\{1,\ldots,N} \sup_{(t_1,t_2)\in\Delta}
w\ssN_i(t_1,t_2)<\infty$,
then a strong uniform law of large numbers
\eqnb
\eqna{PPLATDILLN}
\limf{N}\sup_{(t_1,t_2)\in \Delta} 
\biggl|\frac1N \sum_{i=1}^N
(Z\ssN_i(t_1,t_2)-\EE{Z\ssN_i}(t_1,t_2))\biggr|=0,\ a.e.,
\eqne
holds for the number of arrival times in the intervals
\eqnb
\eqna{PPLATDIincrement}
Z\ssN_i(t_1,t_2):=Z\ssN_i(t_2)-Z\ssN_i(t_1),
\ \ 0\le t_1\le t_2\le T.
\eqne
\DDD
\thme
\prfb
Note that $Z\ssN_i(0)=0$.
It is known \cite[\S 3, Thm.\ 5]{fluid14} that
\[\arrb{l}\dsp
 \EE{Z\ssN_i(t)\, (Z\ssN_i(t)-1)\cdots(Z\ssN_i(t)-q+1)}
\\ \dsp {}
\le\biggl(\int_0^t \max_{s\in[0,u]} w\ssN_i(s,u) \,du\biggr)^{q}
\le (CT)^{q},
\arre\]
holds for all positive integer $q$.
For each positive integers $q$, $N$, and $i$,
$Z\ge 2q$ implies $\dsp Z-1>Z-2>\cdots>Z-q+1>\frac12 Z>0$, so that
\eqnb
\eqna{PPLATDILLNprf1}
\arrb{l}\dsp
\EE{|Z\ssN_i(T)|^{q}}
\\ \dsp {}
=\EE{Z\ssN_i(T)^{q};\ Z\ssN_i(T)\ge 2q}+
\EE{Z\ssN_i(T)^{q};\ Z\ssN_i(T)< 2q}
\\ \dsp {}
\le 2^{q} \EE{Z\ssN_i(T)\cdots(Z\ssN_i(T)-q+1)}+(2q)^{q}
\\ \dsp {}
\le (2CT)^{q}+(2q)^{q}.
\arre
\eqne

Let $0\le s\le t\le T$ and put
$\dsp \Omega(s,t)=\int_s^t w(s,u)\,du$.
According to \cite[\S3]{fluid14}, we have explicit formulas
\[ \begin{array}{l}\dsp
\EE{Z\ssN_i(t)-Z\ssN_i(s)}
\\ \dsp {}
=\sum_{k=1}^{\infty} \int_s^t \biggl( \int_{0\le u_1\le \cdots \le u_k}
w(u_{k-1},u_k)e^{-\Omega(u_{k-1},u_k)}
\\ \dsp \phantom{=\sum\int\int} \times
\prod_{i=1}^{k-1} \left.w(u_{i-1},u_i) e^{-\Omega(u_{i-1},u_i)}
\right|_{u_0=0}\, du_i\biggr)\,du_k
\end{array} \]
and
\[ \begin{array}{l}\dsp
e^{-\Omega(0,t)}
\\ \dsp {}
+\sum_{k=1}^{\infty}
\int_{0\le u_1\le u_2\le \cdots\le u_k\le t} e^{-\Omega(u_k,t)}\,
\prod_{i=1}^k \left. w(u_{i-1},u_i)\,e^{-\Omega(u_{i-1},u_i)}\right|_{u_0=0}
\, du_i
\\ \dsp
=1.
\end{array} \]
Using nonnegativity of each terms and an assumption
$w\ssN_i(u_{k-1},u_k)\le C$
\eqnb
\eqna{PPLATDILLNprf2}
\EE{Z\ssN_i(t)-Z\ssN_i(s)}\le C\int_s^t 1\,du\le C\,(t-s).
\eqne
The estimates \eqnu{PPLATDILLNprf1} and \eqnu{PPLATDILLNprf2},
and \thmu{main2}
in \secu{fSRPvarphiLLNdemonstration} with $M\ssN=1$, $w\ssN_i=C$, and 
$r=1$ imply \eqnu{PPLATDILLN}.
\QED
\prfe


\begin{thebibliography}{99}

\bibitem{Billingsley} 
P.~Billingsley, 
\textit{
Convergence of Probability Measures,
}
2nd ed., John Wiley \& Sons, 1999.

\bibitem{ChowTeicher} 
Y.~S.~Chow, H.~Teicher,
\textit{
Probability theory, independence, interchangeability, martingales,
}
3rd ed., Springer, 2003.

\bibitem{Er49}
P. Erd\"os, \textit{On a Theorem of Hsu and Robbins,}
Ann.\ Math.\ Stat.\ \textbf{20} (1949) 286--291.

\bibitem{Er50}
P.~Erd\"os, \textit{Remark on my Paper `On a Theorem of Hsu and Robbins',}
Ann.\ Math.\ Stat.\ \textbf{21} (1950) 138.

\bibitem{fluid14} 
T.~Hattori,
\textit{
Point process with last-arrival-time dependent intensity and 
1-dimensional incompressible fluid system with evaporation,
}
\url{http://arxiv.org/abs/1409.5117}, Funkcialaj Ekvacioj (2017), to appear.


\bibitem{HR47}
P.~L.~Hsu, H.~Robbins, 
\textit{Complete convergence and the law of large numbers,} 
Proc.\ Nat.\ Acad.\ Sci.\ U.S.A.\ \textbf{33} (1947) 25--31.

\bibitem{TH87}
R.~L.~Taylor, T.~C.~Hu, 
\textit{Strong laws of large numbers for arrays of 
rowwise independent random elements,}
Internat.\ J.\ Math.\ Math. Sci.\ \textbf{10} (1987) 805--814. 

\end{thebibliography}
\end{document}